\documentclass[final,1p,times]{elsarticle}

\usepackage{amssymb}
 \usepackage{amsthm}
\usepackage{amscd}
\usepackage{amsmath}
\usepackage{amsfonts}
\usepackage{amssymb}
\usepackage{graphicx}
\numberwithin{equation}{section}

\newcommand{\ra}{\rightarrow}
\newcommand{\p}{\partial}
\newcommand{\f}{\frac}

\newcommand{\be}{\begin{equation}}
\renewcommand{\ra}{\rightarrow}
\newcommand{\ee}{\end{equation}}
\newcommand{\bea}{\begin{eqnarray}}
\newcommand{\eea}{\end{eqnarray}}
\newcommand{\bna}{\begin{eqnarray*}}
\newcommand{\ena}{\end{eqnarray*}}

\renewcommand{\le}{\left}
\newcommand{\ri}{\right}

\journal{Ann Glob Anal Geom}

\begin{document}

\begin{frontmatter}

\title{Smoothing metrics on closed Riemannian manifolds through the Ricci flow}
\author{Yunyan Yang}
\ead{yunyanyang@ruc.edu.cn}
\address{Department of Mathematics,
Renmin University of China, Beijing 100872, P. R. China}
\begin{abstract}
Under the assumption of the uniform local Sobolev inequality, it is
proved that Riemannian metrics with an absolute Ricci curvature
 bound and a small Riemannian curvature integral bound can be
smoothed to having a sectional curvature bound. This partly extends
previous a priori estimates of Ye Li (J. Geom. Anal. 17 (2007)
495-511; Advances in Mathematics 223 (2010) 1924-1957).
\end{abstract}


\end{frontmatter}

\section{Introduction}
 If a Riemannian manifold has bounded sectional curvature, then its geometric structure is better understood
 than that with weaker curvature bounds, say Ricci curvature bounds.
 Thus it is of significance to deform or smooth a Riemannian metric with a
 Ricci curvature bound to a metric with a sectional curvature bound.
 One way to do this is using the Ricci flow. In this regard we refer the reader to
 the pioneer works \cite{MinO,Dai-Wei-Ye,DYang1,DYang}.
 If the initial
 metric has bounded curvatures, one can show the short time
 existence of the Ricci flow and obtain the covariant derivatives
 bounds for the curvature tensors along the Ricci flow
 \cite{Bando,Shi}. If the initial metric has bounded Ricci
 curvature, under some additional assumption on conjugate radius,
 Dai, etc. studied how to deform the metric on closed manifolds
 \cite{Dai-Wei-Ye}. Also one can deform a metric locally by using
 the local Ricci flow \cite{Li1,Li2,Li3,Guo,DYang}.
 Throughout this paper, we use ${\rm Rm}(g)$ and ${\rm Ric}(g)$
  to denote the Riemannian curvature tensor and Ricci tensor with respect to the metric
  $g$ respectively. Our main result is the following:
  \\

\noindent{\bf Theorem 1.1.} {\it Assume $(M,g_0)$ is a closed
Riemannian manifold of dimension $n$ ($n\geq 3$) and $|{\rm
Ric}(g_0)|\leq K$ for some constant $K$. Let $B_r(x)$ be a geodesic
ball centered at $x\in M$ with radius $r$. Suppose there exists a
constant $A_0>0$ such that for all $x\in M$ and some $r\leq
\min(\f{1}{2}{\rm diam}(g_0),1)$
 \be\label{c1}\le(\int_{B_r(x)}|u|^{\f{2n}{n-2}}dv_{g_0}\ri)^{(n-2)/n}\leq A_0\int_{B_r(x)}|\nabla_{g_0} u|^2dv_{g_0},\quad
 \forall u\in C_0^\infty(B_r(x)).\ee
 Then there exist constants $\epsilon$, $c_1$, $c_2$ depending only on $n$ and $K$ such that if
 \be\label{c2}\le(\int_{B_r(x)}|{\rm
Rm}(g_0)|^{\f{n}{2}}dv_{g_0}\ri)^{{2}/{n}}\leq
{\epsilon}{A_0^{-1}}\quad{\rm for\,\,all}\,\, x\in M,\ee
  then
the Ricci flow
  \be\label{RF}
   \le\{\begin{array}{lll}
   \displaystyle\f{\p g}{\p t}&=-2Ric(g),\\[1.2ex]
   g(0)&=g_0
   \end{array}\ri.
  \ee
  has a unique smooth solution satisfying the following estimates
  \bea\label{g-e}
  |g(t)-g_0|_{g_0}&\leq& c_2t^{\f{2}{n+2}},\\\label{g-e1}
  |{\rm Rm}(g(t))|_\infty&\leq& c_2t^{-1},\\\label{g-e2}
  |{\rm Ric(g(t))}|_\infty&\leq& c_2t^{-\f{n}{n+2}}
  \eea
  for $0\leq t\leq T$ with $T\geq c_1\min(r^2,K^{-1})$.}\\

  When $(M,g_0)$ is a complete noncompact Riemannian manifold,
  similar results were obtained by Ye Li \cite{Li3} and G. Xu \cite{Guo}.
  The assumptions of \cite{Li3}
  is much weaker than (\ref{c1}) and (\ref{c2}) in case $n=4$. It
  comes from Cheeger and Tian's work \cite{Cheeger-Tian}
  concerning the collapsing Einstein 4-manifolds. Here Theorem 1.1 is just the beginning of extending
  the results \cite{Cheeger-Tian,Li2,Li3}, which may depend on the
  Gauss-Bonnet-Chern formula, to general dimensional case.

  For the proof of Theorem 1.1, we follow the lines of \cite{Dai-Wei-Ye,Dai-Wei-Ye1,Li3,DYang}. Let's roughly describe the
  idea. First it is well known \cite{Hamilton,DeTurk} that the Ricci flow (\ref{RF})
  has a unique smooth solution $g(t)$ for a very short time
  interval. Using Moser's iteration and Gromov's covering argument, we derive a priori
  estimates on ${\rm Rm}(g(t))$ and ${\rm Ric}(g(t))$. Let $[0,T_{\rm max}]$ be a maximum time interval
  on which $g(t)$ exists. Then based on those a
  priori estimates, $T_{\rm max}$ has the desired lower bound.

Such kind of results are very useful when considering the relation
between curvature and topology \cite{Anderson,Dai-Wei-Ye,Li2}. Using
Theorem 1.1, we can easily generalize Gromov's almost flat manifold
theorem \cite{Gromov}. Particularly one has the following:\\

\noindent{\bf Theorem 1.2.} {\it There exist constants $\epsilon$
and $\delta$ depending only on $n$ and $K$ such that if a closed
Riemannian manifold $(M,g_0)$ satisfies $|{\rm Ric(g_0)}|\leq K$,
${\rm diam}(g_0)\leq \delta$, (\ref{c1}) and (\ref{c2}) hold for all
$x\in M$, then the universal covering space of $(M,g_0)$ is
$\mathbb{R}^n$. If all the above hypothesis on $(M,g_0)$ are
satisfied and moreover the fundamental
group $\pi(g_0)$ is commutative, then $(M,g_0)$ is diffeomorphic to a torus.}\\

 Before ending this introduction, we would like to mention
\cite{Tian-V} for local regularity estimates for Riemannian
curvatures. The remaining part of the paper is organized as follows.
In Sect. 2, we derive two
 weak maximum principles by using the Moser's iteration. In Sect. 3,
 we estimate the time interval on which the solution of Ricci flow exists, and prove
 Theorem 1.1. Finally Theorem 1.2 is proved in Sect. 4.

\section{Weak maximum principles}
In this section, following the lines of \cite{Li3,DYang}, we give
two maximum principles via the Moser's iteration. Throughout this
section the manifolds need not to be compact. Suppose $(M,g(t))$ are
complete Riemannian manifolds for $0\leq t\leq T$. Let
$\nabla_{g(t)}$ denote the covariant differentiation with respect to
$g(t)$ and $-\Delta_{g(t)}$ be the corresponding Laplace-Beltrami
operator, which will be also denoted by $\nabla$ and $-\Delta$ for
simplicity, the reader can easily recognize it from the context. Let
$A$ be a constant such that for all $t\in[0,T]$, \be\label{Sobolev}
 \le(\int_{B_r(x)}|u|^{\f{2n}{n-2}}dv_{t}\ri)^{(n-2)/n}\leq
 A\int_{B_r(x)}|\nabla u|^2dv_t,\quad \forall u\in
 C_0^\infty(B_r(x)),
\ee where $dv_t=dv_{g(t)}$. Assume that for all $0\leq t\leq T$,
\be\label{g-equ}\f{1}{2}g_0\leq g(t)\leq 2g_0\quad{\rm on}\quad
M.\ee Here and in the sequel, all geodesic balls are defined with
respect to $g_0$.

Firstly we have the following maximum principle:\\

\noindent{\bf Theorem 2.1.} {\it Let $(M,g(t))$ be complete
Riemannian manifolds and (\ref{Sobolev}), (\ref{g-equ}) are
satisfied for $0\leq t\leq T$. Let $f(x,t)$ be such that
\be\label{f-equation}
 \f{\p f}{\p t}\leq \Delta f+uf\quad {\rm on}\quad B_r(x)\times
 [0,T]
\ee with $f\geq 0$, $u\geq 0$, \be\label{dv/dt}
 \f{\p}{\p t}dv_t\leq cu dv_t,
\ee for some constant $c$ depending only on $n$ and for some
$q>n$\be
 \le(\int_{B_r(x)}u^{\f{q}{2}}dv_t\ri)^{\f{2}{q}}\leq \mu
 t^{-\f{q-n}{q}},
\ee where $\mu>0$ is a constant. Then
  for any $p>1$, $t\in[0,T]$, we have
 \be\label{maximum}
  f(x,t)\leq
  CA^{\f{n}{2p}}\le(\f{1+A^{\f{n}{q-n}}\mu^{\f{q}{q-n}}}{t}+\f{1}{r^2}\ri)^{\f{n+2}{2p}}
  \le(\int_0^T\int_{B_r(x)}f^{p}dv_t\ri)^{\f{1}{p}},
 \ee
 where $C$ is a constant depending only on $n$, $q$ and $p$.}\\

 \noindent{\it Proof.} Let $\eta$ be a nonnegative Lipschitz function
 supported in $B_r(x)$. We first consider the case $p\geq 2$. By the
 partial differential inequality (\ref{f-equation}) and (\ref{dv/dt}), we have
 \bna
 \f{1}{p}\f{\p}{\p
 t}\int\eta^2f^pdv_t&\leq&\int\eta^2f^{p-1}\Delta
 fdv_t+C_1
 \int uf^p\eta^2dv_t,
 \ena
 where $C_1$ is a constant depending only on $n$. Integration
 by parts implies
 \bna
 \int\eta^2f^{p-1}\Delta
 fdv_t&=&-2\int\eta f^{p-1}\nabla\eta\nabla fdv_t-(p-1)
 \int\eta^2f^{p-2}|\nabla f|^2dv_t\\
 &=&-\f{4}{p}\int \le(f^{\f{p}{2}}\nabla\eta\nabla(\eta
 f^{\f{p}{2}})-
 |\nabla\eta|^2f^p\ri)dv_t-\f{4(p-1)}{p^2}\\
 &&\times\int\le(|\nabla(\eta f^{\f{p}{2}})|^2+|\nabla\eta|^2f^p
 -2f^{\f{p}{2}}\nabla\eta\nabla(\eta f^{\f{p}{2}})\ri)dv_t\\
 &=&-\f{4(p-1)}{p^2}\int|\nabla(\eta
 f^{\f{p}{2}})|^2dv_t+\f{4}{p^2}\int|\nabla\eta|^2f^pdv_t\\
 &&+\f{4p-8}{p^2}\int f^{\f{p}{2}}\nabla\eta\nabla(\eta
 f^{\f{p}{2}})dv_t\\
 &\leq&-\f{2}{p}\int|\nabla(\eta f^{\f{p}{2}})|^2dv_t+\f{2}{p}
 \int|\nabla\eta|^2f^pdv_t.
 \ena
 Here we have used the elementary inequality $2ab\leq a^2+b^2$. By
 the H\"older inequality, we have
 \bna
  \int uf^p\eta^2dv_t&\leq&\le(\int
  u^{\f{q}{2}}dv_t\ri)^{\f{2}{q}}\le(\int(\eta^2
  f^{p})^{\alpha q_1}dv_t\ri)^{\f{1}{q_1}}\le(\int(\eta^2
  f^{p})^{(1-\alpha) q_2}dv_t\ri)^{\f{1}{q_2}},
 \ena
 where $\f{1}{q_1}+\f{1}{q_2}+\f{2}{q}=1$ and $0<\alpha<1$. Let $\alpha
 q_1=\f{n}{n-2}$ and $(1-\alpha)q_2=1$. This implies
 $q_1=\f{q}{n-2}$, $q_2=\f{q}{q-n}$ and $\alpha=\f{n}{q}$.
 Using the Sobolev inequality (\ref{Sobolev}) and the Young inequality, we obtain
 \bna
  \int uf^p\eta^2dv_t&\leq&\mu t^{-\f{q-n}{q}}\le(\int
  (\eta^2f^p)^{\f{n}{n-2}}dv_t\ri)^{\f{n-2}{q}}\le(\int
  \eta^2f^pdv_t\ri)^{\f{q-n}{q}}\\
  &\leq&\mu t^{-\f{q-n}{q}}\le(A\int
  |\nabla(\eta f^{\f{p}{2}})|^2dv_t\ri)^{\f{n}{q}}\le(\int
  \eta^2f^pdv_t\ri)^{\f{q-n}{q}}\\
  &\leq& \f{1}{pC_1}\int
  |\nabla(\eta
  f^{\f{p}{2}})|^2dv_t+C_2p^{\f{n}{q-n}}\mu^{\f{q}{q-n}}A^{\f{n}{q-n}}t^{-1}\int
  \eta^2f^pdv_t
 \ena
 for some constant $C_2$ depending only on $n$ and $q$.
 Combining all the above estimates one has
 \bea\label{key}
  &&\f{\p}{\p t}\int\eta^2f^pdv_t+\int
  |\nabla(\eta f^{\f{p}{2}})|^2dv_t\leq2\int|\nabla\eta|^2f^pdv_t
  \\&&\quad\quad\quad\quad\quad\quad+C_1C_2p^{\f{q}{q-n}}\mu^{\f{q}{q-n}}A^{\f{n}{q-n}}t^{-1}\int
  \eta^2f^pdv_t.{\nonumber}
 \eea
 For $0<\tau<\tau^\prime<T$, let
 $$\psi(t)=\le\{\begin{array}{lll}
  0, &0\leq t\leq\tau\\[1.2ex]
  \f{t-\tau}{\tau^\prime-\tau},&\tau\leq t\leq \tau^\prime\\[1.2ex]
  1,&\tau^\prime\leq t\leq T.
 \end{array}\ri.$$
 Multiplying (\ref{key}) by $\psi$, we have
  \bea{\nonumber}
  &&\f{\p}{\p t}\le(\psi\int\eta^2f^pdv_t\ri)+\psi\int
  |\nabla(\eta f^{\f{p}{2}})|^2dv_t\leq2\psi\int|\nabla\eta|^2f^pdv_t
  \\&&\quad\quad\quad\quad\quad\quad+\le(C_1C_2p^{\f{q}{q-n}}\mu^{\f{q}{q-n}}A^{\f{n}{q-n}}t^{-1}\psi
  +\psi^\prime\ri)\int\eta^2f^pdv_t.\label{111}
 \eea
 Assume $\tau<\tau^\prime<t\leq T$. Since on the time interval $[\tau,\tau^\prime]$
 $$0\leq\f{\psi(t)}{t}=\f{1}{\tau^\prime-\tau}-\f{\tau}{\tau^\prime-\tau}\f{1}{t}\leq \f{1}{\tau^\prime-\tau}
 \le(1-\f{\tau}{\tau^\prime}\ri)=\f{1}{\tau^\prime},$$
 and on the time interval $[\tau^\prime,T]$ $$\f{1}{T}\leq\f{\psi(t)}{t}\leq
 \f{1}{\tau^\prime},$$
 we have
 \be\label{222}\int_\tau^t\f{\psi(t)}{t}\le(\int \eta^2f^pdv_t\ri)dt\leq \f{1}{\tau^\prime}
 \int_\tau^t\int \eta^2f^pdv_tdt.\ee
 Notice that $0\leq \psi\leq 1$ and $0\leq\psi^\prime\leq
 \f{1}{\tau^\prime-\tau}$.
   Integrating  the differential inequality (\ref{111}) from $\tau$ to
   $t$, we obtain by using (\ref{222})
 \bna
  &&\int\eta^2f^pdv_t+\int_{\tau^\prime}^t\int
  |\nabla(\eta
  f^{\f{p}{2}})|^2dv_tdt\leq2\int_\tau^t\int|\nabla\eta|^2f^pdv_tdt
  \\&&\quad\quad\quad\quad+\le(\f{C_1C_2p^{\f{q}{q-n}}\mu^{\f{q}{q-n}}A^{\f{n}{q-n}}}{\tau^\prime}
  +\f{1}{\tau^\prime-\tau}\ri)\int_\tau^T\int\eta^2f^pdv_tdt.
 \ena
 Applying this estimate and the Sobolev inequality we derive
 \bea\label{itrat}\quad\quad
 \int_{\tau^\prime}^T\int f^{p(1+\f{2}{n})}\eta^{2+\f{1}{n}}dv_tdt&\leq&
 \int_{\tau^\prime}^T\le(\int \eta^2f^pdv_t\ri)^{\f{2}{n}}\le(\int
 f^{\f{pn}{n-2}}\eta^{\f{2n}{n-2}}dv_t\ri)^{\f{n-2}{n}}dt\\{\nonumber}
 &\leq&A\le(\sup_{\tau^\prime\leq t\leq T}\int
 \eta^2f^p\ri)^{\f{2}{n}}\int_{\tau^\prime}^T\int|\nabla(\eta
 f^{\f{p}{2}})|^2dv_tdt\\{\nonumber}
 &\leq&A\le[2\int_\tau^t\int|\nabla\eta|^2f^pdv_tdt+\le(\f{C_1C_2p^{\f{q}{q-n}}
 \mu^{\f{q}{q-n}}A^{\f{n}{q-n}}}{\tau^\prime}\ri.\ri.\\{\nonumber}
 &&\quad\quad\quad\le.\le.+\f{1}{\tau^\prime-\tau}\ri)\int_\tau^T\int\eta^2f^pdv_tdt\ri]^{1+\f{2}{n}}.
 \eea
 For $p\geq p_0\geq 2$ and $0\leq\tau\leq T$, we set
 $$H(p,\tau,r)=\int_\tau^T\int_{B_r(x)}f^pdv_tdt,$$
 where $B_r(x)$ is the geodesic ball centered at $x$ with radius $r$
 measured in $g(0)$. Choosing a suitable cut-off function $\eta$ and
 noting that $|\nabla \eta|_t\leq 2|\nabla\eta|_0$, we obtain from
 (\ref{itrat})
 \bea\label{1}&&
  H\le(p\le( 1+\f{2}{n}\ri),\tau^\prime,r\ri)\\{\nonumber}&&\leq AC_3\le(\f{p^{\f{q}{q-n}}\mu^{\f{q}{q-n}}A^{\f{n}{q-n}}}{\tau^\prime}
  +\f{1}{\tau^\prime-\tau}+\f{1}{(r^\prime-r)^2}\ri)^{1+\f{2}{n}}
  H(p,\tau,r^\prime)^{1+\f{2}{n}},
 \eea
 where $0<r<r^\prime$, $C_3$ is a constant depending only on $n$ and $q$.
 Set
  $$\nu=1+\f{2}{n},\quad p_k=p_0\nu^k,\quad
  \tau_k=(1-\nu^{-\f{qk}{q-n}})t,\quad \quad r_k=(1+\nu^{-\f{qk}{q-n}})r/2.$$
 Then the inequality (\ref{1}) gives
 $$H(p_{k+1},\tau_{k+1},r_{k+1})\leq
 AC_3\le(\f{1+p_0^{\f{q}{q-n}}\mu^{\f{q}{q-n}}A^{\f{n}{q-n}}}{t}+\f{1}{r^2}\ri)^\nu
 \eta^{k\nu}H(p_k,\tau_k,r_k)^\nu,$$
 where $\eta=\nu^{\f{2q}{q-n}}$.
 It follows that
 \bna
  &&H(p_{k+1},\tau_{k+1},r_{k+1})^{\f{1}{p_{k+1}}}\\&&\leq (AC_3)^{\f{1}{p_{k+1}}}
  \le(\f{1+p_0^{\f{q}{q-n}}\mu^{\f{q}{q-n}}A^{\f{n}{q-n}}}{t}+\f{1}{r^2}\ri)^{\f{1}{p_k}}
  \eta^{\f{k}{p_k}}
  H(p_k,\tau_k,r_k)^{\f{1}{p_{k}}}.
 \ena
 Hence we obtain for any fixed $k$
 \bna
 H(p_{k+1},\tau_{k+1},r_{k+1})^{\f{1}{p_{k+1}}}&\leq& (AC_3)^{\sum_{j=0}^k\f{1}{p_{j+1}}}
 \le(\f{1+p_0^{\f{q}{q-n}}\mu^{\f{q}{q-n}}A^{\f{n}{q-n}}}{t}+\f{1}{r^2}\ri)^{\sum_{j=0}^k\f{1}{p_j}}\\
 &&\eta^{\sum_{j=0}^{k}\f{j}{p_j}}H(p_0,\tau_0,r_0)^{\f{1}{p_0}}.
 \ena
 Passing to the limit $k\ra\infty$, one concludes
 \bna
 f(x,t)\leq (CA)^{\f{n}{2p_0}}
 \le(\f{1+(p_0\mu)^{\f{q}{q-n}}A^{\f{n}{q-n}}}{t}+\f{1}{r^2}\ri)^{\f{n+2}{2p_0}}
 \le(\int_0^T\int
 f^{p_0}dv_tdt\ri)^{\f{1}{p_0}}.\ena
 This proves (\ref{maximum}) in the case $p\geq 2$.

 Assuming $f$ satisfies (\ref{f-equation}) and $f\geq 0$. We define a sequence of functions
 $$f_j=f+1/j,\quad j\in\mathbb{N}.$$ Then $f_j$ also satisfies (\ref{f-equation}) and
 $f_j^{p/2}$ is Lipschitz continuous for $1<p<2$. The same argument as
 the case $p\geq 2$ also yields
 \bna
 f_j(x,t)\leq (CA)^{\f{n}{2p_0}}
 \le(\f{1+(p_0\mu)^{\f{q}{q-n}}A^{\f{n}{q-n}}}{t}+\f{1}{r^2}\ri)^{\f{n+2}{2p_0}}
 \le(\int_0^T\int
 f_j^{p_0}dv_tdt\ri)^{\f{1}{p_0}}\ena
 for some constant $C$ depending only on $n$ and $q$, where
 $1<p_0<2$.
 Passing to the limit $j\ra\infty$, we can see that (\ref{maximum})
 holds when $1<p<2$.
 $\hfill\Box$\\

  To proceed we need the following covering lemma belonging to M.
  Gromov.\\

 \noindent{\bf Lemma 2.2} (\cite{Cheeger}, Proposition 3.11). {\it Let $(M,g)$ be a complete Riemannian
 manifold, the Ricci curvature of $M$ satisfy ${\rm Ric}(g)
 \geq (n-1)H$. Then given $r,\epsilon>0$ and $p\in M$, there exists a covering,
 $B_r(p)\subset\cup_{i=1}^NB_\epsilon(p_i)$, ($p_i$ in $B_r(p)$) with $N\leq N_1(n,Hr^2,r/\epsilon)$.
 Moreover, the multiplicity of this covering is at most $N_2(n,Hr^2)$.
 }\\

 \noindent For any complete Riemannian manifold $(M,g_0)$ of dimension $n$
 with
 $|{\rm Ric}(g_0)|\leq K$, it follows from Lemma 2.2 that there exists an absolute
 constant $N$ depending only on $K$ and $n$ such that
 \be
 \label{covering}
 B_{2r}(x)\subset\cup_{i=1}^NB_r(y_i),\quad y_i\in B_{\f{3}{2}r}(x).
 \ee
 Suppose (\ref{Sobolev}) and
(\ref{g-equ}) hold for all $x\in M$ and $0\leq t\leq T$, $g(0)=g_0$.
  Let $f(x,t)$ and $u(x,t)$ be two nonnegative functions satisfying
 $$ \f{\p f}{\p t}\leq \Delta f+C_0f^2,\quad \f{\p u}{\p t}\leq \Delta
 u+C_0fu$$
  on $M\times[0,T]$.
 Assume that there hold on $M\times[0,T]$
 $$u\leq
 c(n)f,\quad\f{\p}{\p t}dv_t\leq c(n)fdv_t.$$
  Define
 \be\label{e00}e_0(t)=\sup_{x\in M,\, 0\leq \tau\leq t}
 \le(\int_{B_{r/2}(x)}f^{\f{n}{2}}dv_{\tau}\ri)^{{2}/{n}}.\ee
 Then we have the following proposition of $f$ and $u$.\\

\noindent {\bf Proposition 2.3.} {\it  Let $f$ and $u$ be as above,
$A$ be given by (\ref{Sobolev}) and $e_0(t)$ be defined by
(\ref{e00}). Suppose there holds for all $x\in M$
 $$\le(\int_{B_{r/2}(x)}f_0^{\f{n}{2}}dv_0\ri)^{\f{2}{n}}\leq (2N^{1+\f{2}{n}}n(C_0+c(n))A)^{-1},$$
 where $N=N(n,K)$ is given by (\ref{covering}), $f_0(x)=f(x,0)$ and $dv_0=dv_{g_0}$. Then there exist two
 constants $C_1$ and $C_2$ depending only on $n$ and $C_0$ such that
 if $0<t<\min(T,C_2N^{-1}r^2)$, then
 $f(x,t)\leq C_1t^{-1}$
 and
 \bna
 u(x,t)\leq
 C_1A^{\f{n}{n+2}}t^{-\f{n}{n+2}}\le[\le(\int_{B_{r}(x)}u_0^{\f{n+2}{2}}dv_0\ri)^{\f{2}{n+2}}+{r^{-\f{4}{n+2}}}
 e_0(t)\ri].
 \ena
  }

 \noindent{\it Proof.} Let $[0,T^\prime]\subset [0,T]$ be the maximal
 interval such that
 \be\label{e0}e_0(T^\prime)=\sup_{x\in M,\, 0\leq t\leq T^\prime}
 \le(\int_{B_{r/2}(x)}f^{\f{n}{2}}dv_t\ri)^{\f{2}{n}}\leq
 ((C_0+c(n))nNA)^{-1}.\ee
 For any cut-off function $\phi$ supported in $B_{r}(x)$, using the same method of
 deriving
 (\ref{key}), we calculate when $p\leq n$ and $m\leq n$,
 \bna
 \f{1}{p}\f{\p}{\p t}\int \phi^{m+2}f^pdv_t
 &\leq& \int\phi^{m+2}f^{p-1}(\Delta
 f+C_0f^2)dv_t+\f{c(n)}{p}\int\phi^{m+2}f^{p+1}dv_t\\
 &\leq& -\int\nabla(\phi^{m+2}f^{p-1})\nabla
 fdv_t+\le(C_0+\f{c(n)}{p}\ri)\\
 &&\quad\times\le(\int_{B_{2r}(x)}f^{\f{n}{2}}dv_t\ri)^{\f{2}{n}}
 \le(\int(\phi^{m+2}f^p)^{\f{n}{n-2}}dv_t\ri)^{\f{n-2}{n}}\\
  &\leq&-\f{2}{p}\int|\nabla(\phi^{\f{m}{2}+1} f^{\f{p}{2}})|^2dv_t+\f{2}{p}
 \int|\nabla\phi^{\f{m}{2}+1}|^2f^pdv_t\\
 &&+\le(C_0+\f{c(n)}{p}\ri)Ne_0A\int|\nabla(\phi^{\f{m}{2}+1}f^{\f{p}{2}})|^2dv_t.\\
 &\leq&-\f{1}{p}\int|\nabla(\phi^{\f{m}{2}+1} f^{\f{p}{2}})|^2dv_t+
 \f{(m+2)^2}{2p}|\nabla\phi|_\infty^2\int\phi^mf^pdv_t.
 \ena
 Here in the second and third inequalities we used (\ref{covering}) and the Sobolev inequality. Hence
 \be\label{2}
 \f{\p}{\p t}\int\phi^{m+2}f^pdv_t+\int|\nabla(\phi^{\f{m}{2}+1}
 f^{\f{p}{2}})|^2dv_t\leq \f{(m+2)^2}{2}|\nabla\phi|_\infty^2\int\phi^mf^pdv_t.
 \ee
 Take $\phi$ supported in $B_{r}(x)$, which is 1 on $B_{r/2}(x)$ and $|\nabla_{g_0} \phi|_\infty^2\leq
 5/r^2$. Since $\f{1}{2}g_{ij}(0)\leq g_{ij}(t)\leq 2g_{ij}(0)$, we have $|\nabla_{g(t)} \phi|_\infty^2\leq
 10/r^2$. Taking $p=\f{n}{2}$ in (\ref{2}) and integrating it from $0$ to $t$, we
 obtain by using (\ref{covering}) again
 \bea
 \int_{B_{r/2}(x)}f^{\f{n}{2}}dv_t{\nonumber}&\leq&
 \int_{B_{r}(x)}f_0^{\f{n}{2}}dv_0+\f{2(m+2)^2}{r^2}\int_0^t\int\phi^mf^{\f{n}{2}}dv_tdt\\
 \label{brx}
  &\leq& N\le(2N^{1+\f{2}{n}}n(C_0+c(n))A\ri)^{-\f{n}{2}}+{2(m+2)^2}{r^{-2}}N(e_0(t))^{\f{n}{2}}t.
 \eea
 Noting that $x$ is arbitrary, one concludes
 $$\le(1-{2(m+2)^2}{r^{-2}}Nt\ri)(e_0(t))^{\f{n}{2}}\leq N\le(2N^{1+\f{2}{n}}n(C_0+c(n))A\ri)^{-\f{n}{2}}.$$
 If $T^\prime<\f{r^2}{8(m+2)^2N}$, then for all $t\in[0,T^\prime]$
 $$e_0(t)<\le(\f{4}{3}\ri)^{{2}/{n}}\le(2Nn(C_0+c(n))A\ri)^{-1}.$$
 This contradicts the maximality of $[0,T^\prime]$. We can therefore
 assume that $T^\prime\geq\min(C_2N^{-1}r^2,T)$.

 It follows from (\ref{2}) that
 \bna
  \f{\p}{\p t}\le(t\int\phi^{m+2}f^pdv_t\ri)&=&t\f{\p}{\p
  t}\int\phi^{m+2}f^pdv_t+\int\phi^{m+2}f^pdv_t\\
  &\leq&\le(\f{(m+2)^2}{2}|\nabla\phi|_\infty^2t+1\ri)\int\phi^mf^pdv_t.
 \ena
 When $0\leq t\leq \min(C_2N^{-1}r^2,T)$,  integrating the above inequality from $0$ to $t$, we have
 \bea\nonumber\int\phi^{m+2}f^pdv_t&\leq&
 \le(\f{2(m+2)^2}{r^2}+\f{1}{t}\ri)\int_0^t\int\phi^mf^pdv_tdt\\\label{3}
 &\leq&c\,t^{-1}\int_0^t\int\phi^mf^pdv_tdt\eea
 for some constant $c$ depending only on $n$.
 Moreover, integrating (\ref{2}) from $0$ to $t$, we derive
 \be\label{4}\int_0^t\int|\nabla (\phi^{\f{m}{2}+1}f^{\f{p}{2}})|^2dv_tdt
 \leq \int\phi^{m+2}f_0^pdv_0+
 \f{2(m+2)^2}{r^2}\int_0^t\int\phi^mf^pdv_tdt.\ee
 Noting that $\f{1}{r^2}\leq \f{C_2}{Nt}$ and $m\leq n$, we calculate by using (\ref{3}) and (\ref{4})
 \bna
 \int_{B_{r/2}(x)}f^{\f{n}{2}+1}dv_t&\leq&\int_{B_{r}(x)}\phi^{m+4}f^{\f{n}{2}+1}dv_t\\
 &\leq&Ct^{-1}\int_0^t\int\phi^{m+2}f^{\f{n}{2}+1}dv_tdt\\
 &\leq&Ct^{-1}\int_0^t\le(\int_{B_{r}(x)}f^{\f{n}{2}}dv_t\ri)^{\f{2}{n}}
 \le(\int(\phi^{m+2}f^{\f{n}{2}})^{\f{n}{n-2}}dv_t\ri)^{\f{n-2}{n}}dt\\
 &\leq&Ct^{-1}N^{\f{2}{n}}e_0(t)A\int_0^t\int|\nabla(\phi^{\f{m}{2}+1}f^{\f{n}{4}})|^2dv_tdt\\
 &\leq&Ct^{-1}N^{\f{2}{n}}e_0(t)A(N(e_0(t))^{\f{n}{2}}+N(e_0(t))^{\f{n}{2}}t)\\
 &\leq& CN^{1+\f{2}{n}}A(e_0(t))^{1+\f{n}{2}}t^{-1},
 \ena
 or equivalently
 \be\label{mut-1}\le(\int_{B_{r/2}(x)}f^{\f{n+2}{2}}dv_t\ri)^{\f{2}{n+2}}\leq
 CNA^{\f{2}{n+2}}e_0(t)t^{-\f{2}{n+2}},\ee
 where $C$ is a constant depending only on $n$, here and in the
 sequel, we often denote various constants by the same $C$.
 Setting $q=n+2$, $p=\f{n}{2}$ and $\mu=CNA^{\f{2}{n+2}}e_0(T^\prime)$, we obtain by employing Theorem 2.1
  \bna
  f(x,t)&\leq& CA\le(\f{1+A^{\f{n}{2}}\mu^{\f{n+2}{2}}}{t}+\f{1}{r^2}\ri)^{\f{n+2}{n}}
  \le(\int_0^{t}\int_{B_r(x)}f^{\f{n}{2}}dv_tdt\ri)^{\f{2}{n}}\\
  &\leq&CAe_0(T^\prime){t}^{\f{2}{n}}
  \le(\f{1+A^{\f{n}{2}}\mu^{\f{n+2}{2}}}{t}+\f{1}{r^2}\ri)^{\f{n+2}{n}}
  \ena
  for $t\in [0,T^\prime]$.
  Recalling the definition of $e_0(T^\prime)$ (see (\ref{e0}) above), we can
  see that $Ae_0(T^\prime)$ is bounded and
 \be\label{mu}A^{\f{n}{2}}\mu^{\f{n+2}{2}}=(CNAe_0(T^\prime))^{\f{n+2}{2}}\ee
 is also bounded. Therefore when $0<t<\min(T,C_2N^{-1}r^2)$, $f(x,t)\leq C_1t^{-1}$
 for some constants $C_1$ and $C_2$ depending only on $n$, $C_0$.

 Using $u\leq c(n)f$ and $\partial_tdv_t\leq c(n)fdv_t$ and mimicking the method of proving (\ref{2}), we obtain
 \be\label{5}
 \f{\p}{\p t}\int\phi^{m+2}u^pdv_t+\int|\nabla(\phi^{\f{m}{2}+1}
 u^{\f{p}{2}})|^2dv_t\leq \f{C}{r^2}\int\phi^mu^pdv_t.
 \ee
  Taking $m=0$, $p=n/2$ and integrating this inequality, we have by
 using (\ref{covering}))
 \be\label{6}\int_0^t\int|\nabla(\phi u^{\f{n}{4}})|^2dv_tdt\leq \int_{B_{r}(x)}u_0^{\f{n}{2}}dv_0
 +\f{C}{r^2}N(e_0(t))^{\f{n}{2}}t.\ee
Integrating (\ref{5}) with $m=2$, $p={(n+2)}/{2}$, and using the
Sobolev inequality
 (\ref{Sobolev}), we obtain
 \bna
 \int_{B_{r/2}(x)}u^{\f{n+2}{2}}dv_t&\leq&
 \int_{B_{r}(x)}u_0^{\f{n+2}{2}}dv_0+\f{C}{r^2}\int_0^t\int\phi^2u^{\f{n+2}{2}}dv_tdt\\
 &\leq&\int_{B_{r}(x)}u_0^{\f{n+2}{2}}dv_0+\f{C}{r^2}e_0(t)A\int_0^t\int|\nabla(\phi
 u^{\f{n}{4}})|^2dv_tdt,
 \ena
 which together with (\ref{6}) and (\ref{covering}) gives
 \bea\nonumber \int_{B_{r/2}(x)}u^{\f{n+2}{2}}dv_t&\leq& \int_{B_{r}(x)}u_0^{\f{n+2}{2}}dv_0+\f{C}{r^2}
 e_0(t)A\le(\int_{B_{r}(x)}u_0^{\f{n}{2}}dv_0+\f{C}{r^2}Ne_0(t)^{\f{n}{2}}t\ri)\\
 &\leq&\int_{B_{r}(x)}u_0^{\f{n+2}{2}}dv_0+
 \f{C}{r^2}NA(e_0(t))^{1+\f{n}{2}}\le(1+\f{1}{r^2}t\ri).\label{kk}\eea
 Notice that when $0\leq t\leq \min(C_2r^2/N,T)$, (\ref{mut-1})
 implies
 $$\int_{B_{r/2}(x)}f^{\f{n+2}{2}}dv_t\leq \mu t^{-1}.$$
 Without loss of generality we can assume $A>1$ (otherwise we can substitute $A$ for $A+1$).
 In view of (\ref{mu}) and (\ref{kk}), we obtain by using Theorem 2.1 in the case $q=n+2$ and $p=(n+2)/2$
 \bna
 u(x,t)&\leq&CA^{\f{n}{n+2}}\le(\f{1}{t}+\f{1}{r^2}\ri)
 \le(\int_0^t\int_{B_{r/2}(x)}u^{\f{n+2}{2}}dv_tdt\ri)^{\f{2}{n+2}}\\
 &\leq&CA^{\f{n}{n+2}}t^{-\f{n}{n+2}}\le[\le(\int_{B_{r}(x)}u_0^{\f{n+2}{2}}dv_0\ri)^{\f{2}{n+2}}+{r^{-\f{4}{n+2}}}
 e_0(t)\ri],
 \ena
 provided that $0\leq t\leq \min(C_2r^2/N,T)$.
  $\hfill\Box$\\

  \noindent{\bf Remark 2.4.} We remark that Theorem 2.1 and Proposition 2.3 are very similar to Theorem
  A.1 and Corollary A.10 of Dean Yang's paper \cite{DYang1}
  respectively. The differences are that we have heat flow type
  inequalities, but Dean Yang has heat flow type inequalities with
  cut-off function. It seems that Dean Yang's Corollary A.10 is stronger
  than our Proposition 2.3, which is enough for our use here. Also
  we should compare Theorem 2.1 with (\cite{Dai-Wei-Ye,Dai-Wei-Ye1}, Theorem 2.1), where Dai-Wei-Ye obtained a
  similar result by using a similar method. Here the constant $C$ of
  (\ref{maximum})
  depends only on $n$, $q$, $p$, but not on the Sobolev constant $A$.
  While in \cite{Dai-Wei-Ye,Dai-Wei-Ye1}, since the Sobolev constants $C_S(t)$ along the flow
  are bounded, they need not care how the constant $C$ exactly depends on $C_S$.
 \section{Short time existence of the Ricci flow }
 In this section we focus on closed Riemannian manifolds. Precisely, following the lines
 of \cite{Li3,DYang}, we study the short time existence of the
 Ricci flow and give the proof of Theorem 1.1.  Assume $(M,g_0)$ is a closed Riemannian manifold  of dimension
 $n (\geq 3)$
 with $|{\rm Ric}(g_0)|\leq K$. Consider the Ricci flow
 \be\label{ricciflow}
   \le\{\begin{array}{lll}
   \displaystyle\f{\p g}{\p t}&=-2Ric(g),\\[1.2ex]
   g(0)&=g_0.
   \end{array}\ri.
  \ee
 It is well known \cite{Hamilton} that the Riemannian curvature tensor
 and the Ricci curvature tensor satisfy the following evolution
 equations
  \bea
  \f{\p {\rm Rm}}{\p t}&=&\Delta {\rm Rm}+{\rm Rm}*{\rm Rm},\\
  \f{\p{\rm Ric}}{\p t}&=&\Delta{\rm Ric}+{\rm Rm}*{\rm Ric},
  \eea
  where ${\rm Rm}*{\rm Rm}$ is a tensor that is quadratic in ${\rm
  Rm}$, ${\rm Ric}*{\rm Rm}$ can be understood in a similar way. It
  follows that
  \bea\label{cur}
  \f{\p |{\rm Rm}|}{\p t}&\leq&\Delta|{\rm Rm}|+c(n)|{\rm
  Rm}|^2,\\{\label{Ricur}}
  \f{\p |{\rm Ric}|}{\p t}&\leq&\Delta|{\rm Ric}|+c(n)|{\rm Rm}||{\rm
  Ric}|.
  \eea
 To prove Theorem 1.1, it suffices to prove the following:\\

  \noindent{\bf Proposition 3.1.} {\it Let $(M,g_0)$ be a closed Riemannian manifold of
  dimension $n (\geq 3)$ with  $|{\rm Ric}(g_0)|\leq K$. Suppose there exists a constant
  $A_0>0$ such that the following local Sobolev inequalities
  hold for all $x\in M$
  $$\|u\|_{{2n}/{(n-2)}}^2\leq A_0\|\nabla u\|_2^2,\quad \forall u\in C_0^\infty(B_r(x)).$$
 Then there exist constants $C_1$, $C_3$
depending only on $n$ and $K$, and $C_2$ depending only on $n$ such
that for $r\leq 1$, if
$$\le(\int_{B_{r/2}(x)}|{\rm
Rm}(g_0)|^{\f{n}{2}}dv_{g_0}\ri)^{2/n}\leq (C_1A_0)^{-1}
$$ for all $x\in M$, then the Ricci flow (\ref{ricciflow})
  has a smooth solution for $0\leq t\leq T$, where
  $T\geq C_2\min(r^2/N,K^{-1})$, such that for all $x\in M$
  \bea
  \label{g-equiv}
   &&\f{1}{2}g_0\,\leq\, g(t)\,\leq 2g_0,
  \\\label{R1}
  &&\|u\|_{{2n}/{(n-2)}}^2\leq 4A_0\|\nabla u\|_2^2,\quad \forall u\in C_0^\infty(B_r(x)),\\\label{Ric1}
  &&\le(\int_{B_{r/2}(x)}|{\rm Rm}(g(t))|^{\f{n}{2}}dv_t\ri)^{2/n}\leq 2N(C_1A_0)^{-1}.
  \eea}

  \noindent{\it Proof.} It is
 well known (see for example \cite{DeTurk, Hamilton}) that a smooth solution $g(t)$ of the Ricci
 flow (\ref{ricciflow})
exists for a short time interval and is unique.
  Let $[0,T_{\rm max})$ be a maximum time interval on which
  $g(t)$ exists and (\ref{g-equiv})-(\ref{Ric1}) hold. Clearly
  $T_{\rm max}>0$ since the strict inequalities in
  (\ref{g-equiv})-(\ref{Ric1}) hold at $t=0$.
  Suppose $T_{\rm max}<T_0=C_2\min(r^2/N,K^{-1})$ for some constant
  $C_2$ to be determined later.
  Since the Ricci curvature satisfies (\ref{Ricur}), it follows from
Proposition 2.3 that for $0\leq t\leq T^\prime$,
  \bea{\nonumber}
  |{\rm Ric}(g(t))|&\leq&CA_0^{\f{n}{n+2}}t^{-\f{n}{n+2}}\le[\le(\int_{B_r(x)}
  |{\rm
  Ric}(g_0)|^{\f{n+2}{2}}dv_0\ri)^{\f{2}{n+2}}
  +{r^{-\f{4}{n+2}}}e_0(T^\prime)\ri]\\{\nonumber}
  &\leq&CA_0^{\f{n}{n+2}}t^{-\f{n}{n+2}}\le(K^{\f{2}{n+2}}(e_0(T^\prime))^{\f{n}{n+2}}+{r^{-\f{4}{n+2}}}
  e_0(T^\prime)\ri)\\\label{RicBound}
  &\leq&C(K^{\f{2}{n+2}}+r^{-\f{4}{n+2}})t^{-\f{n}{n+2}},
  \eea
  where $T^\prime$ and $e_0(T^\prime)$ are defined by
  (\ref{e0}) in the case $f$ is replaced by $|{\rm Rm}|$.
  It follows that for all $x\in M$, $u\in C_0^\infty(B_r(x))$ and $0\leq t\leq
  T^\prime$,
  \bna
  \le|\f{d}{dt}\int_{B_r(x)} |u|^{\f{2n}{n-2}}dv_t\ri|&\leq& 2|{\rm Ric}(g(t))|_\infty\int_{B_r(x)}
  |u|^{\f{2n}{n-2}}dv_t\\
  &\leq& Ct^{-\f{n}{n+2}}\int_{B_r(x)} |u|^{\f{2n}{n-2}}dv_t.
  \ena
  This implies
  $$e^{-C\,t^{\f{2}{n+2}}}\int_{B_r(x)} |u|^{\f{2n}{n-2}}dv_0\leq \int_{B_r(x)} |u|^{\f{2n}{n-2}}dv_t
  \leq e^{C\,t^{\f{2}{n+2}}}\int_{B_r(x)} |u|^{\f{2n}{n-2}}dv_0.$$
  Similarly we have
  $$\le|\f{d}{dt}\int_{B_r(x)} |\nabla u|^{2}dv_t\ri|\leq Ct^{-\f{n}{n+2}}\int_{B_r(x)} |\nabla u|^{2}dv_t,$$
  and
  $$e^{-C\,t^{\f{2}{n+2}}}\int_{B_r(x)} |\nabla u|^{2}dv_0\leq \int_{B_r(x)} |\nabla u|^{2}dv_t
  \leq e^{C\,t^{\f{2}{n+2}}}\int_{B_r(x)} |\nabla u|^{2}dv_0.$$
  Hence if $T_{\rm max}< T_0=C_2\min(r^2/N,K^{-1})$ for sufficiently small $C_2$ depending only on
  $n$ and $K$, then (\ref{R1}) holds with strict inequality.

  To show (\ref{g-equiv}) holds with strict inequality, we fix a tangent vector $v$ and calculate
  \bna
  \f{d}{dt}|v|_{g(t)}^2=\f{d}{dt}(g_{ij}(t)v^iv^j)=-2{\rm
  Ric}_{ij}v^iv^j,
  \ena
  which together with (\ref{RicBound}) gives
  $$\le|\f{d}{dt}\log|v|_{g(t)}^2\ri|\leq C(K^{\f{2}{n+2}}+r^{-\f{4}{n+2}})t^{-\f{n}{n+2}}.$$
  Therefore we obtain for $0\leq t<C_2\min(r^2,K^{-1})$,
  $$\f{1}{2}|v|_{g(0)}^2<|v|_{g(t)}^2<2|v|_{g(0)}^2.$$
  Using the same method of deriving (\ref{brx}), one can see that
  the strict inequality in (\ref{Ric1}) holds
  when $0\leq t<C_2\min(r^2,K^{-1})$ for sufficiently small $C_2$.
  By Proposition 2.3, $|{\rm Rm}(g(t))|_\infty\leq Ct^{-1}$ for all $t\in [0,T_{\rm max}]$. Hence
 one can extend $g(t)$ smoothly beyond $T_{\rm max}$ with
 (\ref{g-equiv})-(\ref{Ric1})
 still holding. This contradicts the assumed
maximality of $T_{\rm max}$. Therefore $T_{\rm max}\geq T_0$.
$\hfill\Box$\\

\noindent{\it Proof of Theorem 1.1.} By Proposition 3.1, there
exists a unique solution $g(t)$ of the Ricci flow (\ref{ricciflow})
such that (\ref{g-equiv})-(\ref{Ric1}) hold. Then by Proposition
2.3, one concludes
 $$|{\rm Rm}(g(t))|\leq Ct^{-1},\quad |{\rm Ric}(g(t))|\leq Ct^{-\f{n}{n+2}}$$
 for $t\in [0,T_0]$. This completes the proof of Theorem 1.1.
 $\hfill\Box$

\section{Applications}

In this section, we will prove Theorem 1.2 by applying Theorem 1.1.
It follows from (\ref{g-e})-(\ref{g-e2}) that the deformed metric
$g(t)$ has uniform sectional curvature bounds away from $t=0$ and
$g(t)$ is close to $g(0)$ when $t$ is close to $0$. We
first show that diameters of the flow are under control, namely\\

\noindent{\bf Lemma 4.1.} {\it Let $g(t)$ be the Ricci flow in
Theorem 1.1. Then for $0\leq t\leq c_1\min(r^2,K^{-1})$, there
exists a constant $c$ depending only on $n$ and $K$ such that
 \be \label{diam}
  e^{-ct^{\f{2}{n+2}}}{\rm diam}(g_0)\leq {\rm diam}(g(t))\leq
  e^{ct^{\f{2}{n+2}}}{\rm diam}(g_0).
 \ee
where ${\rm diam}(g(t))$ means the diameter of the manifold $(M,g(t))$.}\\

\noindent{\it Proof.} Let $\gamma:[0,1]\ra M$ be any smooth curve.
Denote the length of $\gamma$ by
$$l_\gamma(t)=\int_0^1|\dot{\gamma}(s)|_{g(t)}^2ds.$$
We calculate by using the Ricci bound in Theorem 1.2
 \bna
 \le|\f{d}{dt}l_\gamma(t)\ri|= \le|
 \int_0^1-2{\rm Ric}_{g(t)}(\dot{\gamma}(s),\dot{\gamma}(s))ds\ri|
\leq ct^{-\f{n}{n+2}}l_\gamma(t). \ena This implies
$$l_\gamma(0)e^{-ct^{\f{2}{n+2}}}\leq l_\gamma(t)\leq l_\gamma(0)e^{ct^{\f{2}{n+2}}}.$$
It follows that
$$e^{-ct^{\f{2}{n+2}}}{\rm dist}_{g_0}(p,q)\leq {\rm dist}_{g(t)}(p,q)
\leq e^{ct^{\f{2}{n+2}}}{\rm dist}_{g_0}(p,q),$$ where ${\rm
dist}_{g(t)}(p,q)$ denote the distance between $p$ and $q$ in the
metric $g(t)$. This gives the desired result. $\hfill\Box$\\

The following proposition is a corollary of Gromov's almost flat manifold theorem \cite{Gromov}:\\

\noindent{\bf Proposition 4.2 (Gromov).} {\it Let $(M,g)$ be a
compact Riemannian manifold of dimension $n$. Assume the sectional
curvature is bounded, i.e., $|{\rm Sec}(g)|\leq \Lambda$. Then there
exists a constant $\epsilon_0$ depending only on $n$ such that if
\be\label{flat}\Lambda ({\rm diam}(g))^2\leq \epsilon_0,\ee then the
universal covering of $(M,g)$ is diffeomorphic to $\mathbb{R}^n$. If
in addition the fundamental group $\pi(M)$ is
commutative, then $(M,g)$ is diffeomorphic to a torus.}\\

\noindent{\it Proof of Theorem 1.2.} Let $g(t)$ be a unique solution
to the Ricci flow (\ref{RF}). By (\ref{g-e1}), for $0\leq t\leq
c_1\min(r^2,K^{-1})$,
$$|{\rm Sec}(g(t))|\leq ct^{-1},$$
 where ${\rm Sec}(g(t))$ denotes the sectional curvature of $(M,g(t))$. Let
 $\epsilon_0$ be given by Proposition 4.2. Take $t_0=c_1\min(r^2,K^{-1})$ and
 $$\delta=\le(\epsilon_0t_0c^{-1}e^{-2ct_0^{\f{2}{n+2}}}\ri)^{1/2}.$$
 If ${\rm diam}(g_0)\leq \delta$, then we obtain by Lemma 4.1
 $$|{\rm Sec}(g(t_0))({\rm diam}(g(t_0)))^2|\leq ct_0^{-1}e^{2ct_0^{\f{2}{n+2}}}
 ({\rm diam}(g_0))^2\leq \epsilon_0.$$
 Applying Proposition 4.2 to $g(t_0)$, we conclude Theorem 1.2.
 $\hfill\Box$\\

 {\bf Acknowledgements.} The author is partly supported by the program for
 NCET. He thanks Ye Li for introducing this interesting topic to him.
 Also he thanks the referee for valuable comments and suggestions,
 which improve this paper.

\bigskip

\end{document}